\documentclass[12pt]{amsart}

\usepackage{fullpage} 
\usepackage{amsmath}
\usepackage{amssymb}
\usepackage{amsthm}

\newtheorem{thm}{Theorem}[section]
\newtheorem{lem}[thm]{Lemma}

\newtheorem{prop}[thm]{Proposition}

\theoremstyle{definition}

\newtheorem{defn}{Definition}

\newcommand{\cat}{{^\smallfrown}}
\newcommand{\res}{\upharpoonright}

\newcommand{\pp}{\mathcal{U}}
\newcommand{\pq}{\mathcal{V}}
\newcommand{\pie}{\Pi_1^0}
\newcommand{\mm}{\mathcal{M}}
\newcommand{\nn}{\mathcal{N}}

\begin{document}

\title{Embedding $FD(\omega)$ into $\mathcal{P}_s$ Densely}

\author{Joshua A. Cole \\University of Notre Dame \\ jcole1@nd.edu}

\thanks{Thanks to my adviser Peter Cholak for his guidance in my research. My research was partially supported by NSF grants DMS-0245167 and RTG-0353748 and a Schmitt Fellowship at the University of Notre Dame.}

\begin{abstract}

Let $\mathcal{P}_s$ be the lattice of degrees of non-empty $\Pi_1^0$ subsets of $2^\omega$ under Medvedev reducibility. Binns and Simpson proved that $FD(\omega)$, the free distributive lattice on countably many generators, is lattice-embeddable below any non-zero element in $\mathcal{P}_s$. Cenzer and Hinman proved that $\mathcal{P}_s$ is dense, by adapting the Sacks Preservation and Sacks Coding Strategies used in the proof of the density of the c.e.\ Turing  degrees. With a construction that is a modification of the one by Cenzer and Hinman, we improve on the result of Binns and Simpson by showing that for any $\mathcal{U} <_s \mathcal{V}$, we can lattice embed $FD(\omega)$ into $\mathcal{P}_s$ strictly between $deg_s(\mathcal{U})$ and $deg_s(\mathcal{V)}$. We also note that, in contrast to the infinite injury in the proof of the Sacks Density Theorem, in our proof all injury is finite, and that this is also true for the proof of Cenzer and Hinman, if a straightforward simplification is made. 

\end{abstract}

\maketitle

\section{Introduction and Basic Definitions}
Recently there has been renewed interest in the study of mass problems, which is the study of the lattice of equivalence classes of subsets of $\omega^\omega$ under either Medvedev or Muchnik reducibility. A subset of $\omega^\omega$ can be thought of as a mathematical problem, and then an element of the subset is thought of as a solution of the problem.  In this paper we focus on Medvedev reducibility: if $\pp$ and $\pq$ are subsets of $\omega^\omega$, then $\pp$ is \emph{Medvedev} (or \emph{strongly}) reducible to $\pq$ $(\pp \leq_s \pq)$ if there is a Turing functional $\Psi$ so that for all $f \in \pq$, $\Psi(f) \in \pp$. The intuition is that if $\pp \leq_s \pq$ then a solution to the mass problem $\pq$ uniformly yields a solution to the mass problem $\pp$. Muchnik reducibility relaxes the uniformity constraint: if $\pp$ and $\pq$ are subsets of $\omega^\omega$, then $\pp$ is \emph{Muchnik} (or \emph{weakly}) reducible to $\pq$ $(\pp \leq_w \pq)$ if for every $f \in \pq$ there is a Turing functional $\Psi$ so that $\Psi(f) \in \pp$. Under either reducibility, there is an equivalence relation defined in the usual way:  $\mathcal{U} \equiv \mathcal{V}$ just in case $\mathcal{U} \leq \mathcal{V}$ and $\mathcal{U} \leq \mathcal{V}$. We write $deg_s(\pp)$ $(deg_w(\pp))$ for the degree of $\pp$ under the equivalence relation induced by strong (weak) reducibility.

\indent
In this paper we work mostly with subsets of $2^\omega$, and this allows for the use of compactness; we also work mostly with $\Pi_1^0$ mass problems, and this allows for the use of computable approximations. Often we will use `$\Pi_1^0$ class' as shorthand for `$\Pi_1^0$ subset of $2^\omega$'. The lattices of non-empty $\Pi_1^0$ subsets of $2^\omega$ under the Medvedev and Muchnik reducibilities will be denoted by $\mathcal{P}_s$ and $\mathcal{P}_w$, respectively. The subscripts `s' and `w' stand for `strong' and `weak', respectively; their use avoids confusion arising from the names of both reducibilities beginning with an `M'. Much research in mass problems focuses specifically on $\mathcal{P}_s$ or $\mathcal{P}_w$.

\indent
In the lattice $\mathcal{P}_s$ the join is induced by an operation on representative $\Pi_1^0$ classes. If $\mathcal{U}$ and $\mathcal{V}$ are $\Pi_1^0$ classes, then \[ \mathcal{U} \vee \mathcal{V} = \lbrace f \oplus g : f \in \mathcal{U}, g \in \mathcal{V} \rbrace. \] 

The meet is also induced by an operation on representative $\Pi_1^0$ classes. It is similar to a union, but we also need to be able to tell which mass problem an element originally came from.  If $\mathcal{U}$ and $\mathcal{V}$ are $\Pi_1^0$ classes, then \[ \mathcal{U} \wedge \mathcal{V} = \lbrace 0 \cat f : f \in \mathcal{U} \rbrace \cup \lbrace 1 \cat g : g \in \mathcal{V} \rbrace,\] where for a function $h$ and number $n$, $n \cat h$ is the function given by $n\cat h(0) = n$ and 
\newline
$n \cat h (x) = h(x-1)$ if $x \neq 0$. 

It is a straightforward exercise to check that the induced operations on equivalence classes are well-defined and that  $\mathcal{P}_s$ together with these induced operations is a lattice and is distributive. As usual when working with a collection of equivalence classes, we will work with representative elements, as we did above in the definitions of the meet and the join. 

\indent
There is a bottom element $\textbf{0}$: it consists of the $\Pi_1^0$ subsets of $2^\omega$ which contain a computable function. There is also a top element $\textbf{1}$. 
\newline
\indent
Why do we study the lattice $\mathcal{P}_s$? There are interesting connections between $\mathcal{P}_s$ and classical mathematical logic. For instance, the set of completions of Peano Arithmetic has degree $\textbf{1}$ in $\mathcal{P}_s$. Perhaps more importantly, $\mathcal{P}_s$ is a refinement of $\mathcal{P}_w$; Simpson and Slaman ~\cite{SiSl} have shown that every non-zero $\mathcal{P}_w$-degree contains infinitely many $\mathcal{P}_s$-degrees. With this  basic framework in mind, it is reasonable to suppose that the study of the local properties of $\mathcal{P}_s$ may shed light on the the local properties of $\mathcal{P}_w$. In turn there are many connections in $\mathcal{P}_w$ to topics in or related to computability theory, such as the c.e.\ Turing degrees ($\mathcal{R}_T$) ~\cite{Si3}, almost everywhere domination ~\cite{Si2}, the diagonally non-recursive functions, randomness, and computational complexity ~\cite{Si1}. The theory of mass problems provides an useful context in which to think about problems in these areas. 
\newline
\indent
A way to investigate the local properties of $\mathcal{P}_w$ and $\mathcal{P}_s$ is to ask the questions that were answered in the case of $\mathcal{R}_T$. An easy first result is that while $\mathcal{R}_T$ is only an upper semi-lattice, $\mathcal{P}_w$ and $\mathcal{P}_s$ form true lattices. Binns ~\cite{Bi} has shown every non-trivial degree splits in both $\mathcal{P}_w$ and $\mathcal{P}_s$, as Sacks showed of $\mathcal{R}_T$.  Alfeld ~\cite{Al} has studied the analogous question in the upwards direction, namely which degrees branch. Sacks proved the density of $\mathcal{R}_T$; Cenzer and Hinman ~\cite{CeHi} proved the density of $\mathcal{P}_s$, but whether $\mathcal{P}_w$ is dense is not known. 
\newline
\indent
Binns and Simpson have studied which lattices embed in $\mathcal{P}_s$ and $\mathcal{P}_w$, and our theorem is an improvement on two of their results in $\mathcal{P}_s$. Binns ~\cite{Bi} proved that every finite distributive lattice embeds densely in $\mathcal{P}_s$. Together Binns and Simpson ~\cite{BiSi} proved that there is a lattice embedding of $FD(\omega)$, the free distributive lattice on countably many generators, below any non-trivial $deg_s(\pq) \in \mathcal{P}_s$. Our result makes this embedding dense.

\begin{thm} (\textbf{Main Result}) If $deg_s(\pp) <_s \deg_s(\pq)$ in $\mathcal{P}_s$, the lattice of degrees of non-empty $\Pi_1^0$ subsets of $2^\omega$ under Medvedev reducibility, then there is a lattice-embedding of $FD(\omega)$, the countable free distributive lattice, strictly between $deg_s(\pp)$ and $deg_s(\pq)$.
\end{thm}

The method of construction for our extension of the results of Binns and Simpson is derived from the proof of the density of $\mathcal{P}_s$ by Cenzer and Hinman. They used separating classes of c.e.\ sets to construct $\Pi_1^0$ classes. They satisfied requirements with Sacks Coding and Preservation Strategies, the main techniques of Sacks' proof of the density of $\mathcal{R}_T$. 

Sacks faced and solved the problem of infinite injury in his proof. Once the proper definitions are made and preliminary lemmas proved, the proof of Cenzer and Hinman follows closely the proof of the Sacks Density Theorem in a style as in, for example, the proof given by Soare ~\cite[142-145]{So}. Like the proof of the density of $\mathcal{R}_T$, the proof of the density of $\mathcal{P}_s$ has infinite injury. Although the construction in this paper is based on the construction of Cenzer and Hinman, with one modification we are able to eliminate the possibility of infinite injury. In fact, the same modification eliminates infinite injury in the construction of Cenzer of Hinman. We will give more details in Sections ~\ref{sec:restraints} and ~\ref{sec:negative}.

Another aspect of our dense-embedding result is that it suggests, along with other evidence, that we can often do anything we want in $\mathcal{P}_s$ densely, if we can do it at all. Simpson and Binns ~\cite{BiSi} showed $FD(\omega)$ embeds in $\mathcal{P}_w$, and then they showed it embeds in $\mathcal{P}_s$. From this point, we make use of known techniques and a finite-injury priority argument to make an embedding of $FD(\omega)$ in $\mathcal{P}_s$ that is dense. 

Binns' paper ~\cite{Bi} on splitting in $\mathcal{P}_w$ and $\mathcal{P}_s$ is another example of this process. First splitting is shown in $\mathcal{P}_w$, and in $\mathcal{P}_s$, and then it is shown to occur densely in $\mathcal{P}_s$, still with only a finite-injury priority argument. Moreover, an attempt to use Binns' methods for the dense splitting in $\mathcal{P}_s$ does not directly yield a dense splitting in $\mathcal{P}_w$. Similarly, the proof from this paper cannot easily be modified to embed $FD(\omega)$ densely in $\mathcal{P}_w$. This is because the length of agreement function used in this paper has no easy, well-behaved analogue in the case of weak reducibility. 

The existence of a dense splitting in $\mathcal{P}_w$, or of a dense embedding of $FD(\omega)$ in $\mathcal{P}_w$, would immediately imply the density of $\mathcal{P}_w$. Thus the result of this paper is a little bit more evidence that new techniques may be needed to answer the question of density for $\mathcal{P}_w$.

\section{The General Plan}
Let $\mathcal{U} <_s \mathcal{V}$. (From here on, we will suppress the subscript $s$ on $<, \leq, \equiv,$ etc.) With a priority argument we will construct sequences of c.e.\ subsets of $\omega$, $\lbrace A_i \rbrace_{i \in \omega}$ and $\lbrace B_i \rbrace_{i \in \omega}$, satisfying certain properties and such that for each $i$, $A_i \cap B_i = \emptyset$. Then for each $i$ we set 
\[ \mathcal{S}_i = S(A_i, B_i) = \lbrace X \in 2^\omega :  n \in A_i \Rightarrow X(n) = 1, n \in B_i \Rightarrow X(n) = 0 \rbrace. \] In general if $A$ and $B$ are disjoint c.e. subsets of $\omega$, then $S(A,B)$ is called a \emph{separating class}. It is a straightforward exercise to show that every separating class is a $\Pi_1^0$ class. Finally, for each $i$ we set \[ \mathcal{V}_i = (\mathcal{U} \vee \mathcal{S}_i) \wedge \mathcal{V}. \]
\newline
\indent
We will consider the lattice $\mathcal{L}$ in $\mathcal{P}_s$ generated by $\lbrace deg_s(\mathcal{V}_i) \rbrace_{i \in \omega}$, and show that if certain requirements are satisfied then $\mathcal{L}$ is free and entirely between $deg_s(\mathcal{U})$ and $deg_s(\mathcal{V})$. Note that the free distributive lattice on countably many generators has no maximal or minimal element.

\section{Requirements}

In our priority construction we will have positive and negative requirements. 

\bigskip

For each pair $I,J$ of finite subsets of $\omega$ such that $I \cap J = \emptyset$ we have the following requirements:
\newline
\[P_{I,J}: \mathcal{U} \vee \bigvee_{i \in I} \mathcal{S}_i \ngeq (\mathcal{U} \vee \bigwedge_{j \in J} \mathcal{S}_j) \wedge \mathcal{V} \] 

\[N_I: \mathcal{U} \vee \bigvee_{i \in I} \mathcal{S}_i \ngeq \mathcal{V}. \] 

\bigskip

We must verify five facts to show that if we satisfy $P_{I,J}$ and $N_I$, then the lattice $\mathcal{L}$ generated by $\lbrace deg_s(\mathcal{V}_i) \rbrace_{i \in \omega}$ is free and between $deg_s(\mathcal{U})$ and $deg_s(\mathcal{V})$. From here on, we assume that $I$ and $J$ are finite subsets of the natural numbers.

\begin{enumerate}

\medskip
\item Every element of $\mathcal{L}$ is below $\deg_s(\mathcal{V})$. 

\medskip
\item Every element of $\mathcal{L}$ is above $\deg_s(\mathcal{U}).$ 

\medskip
\noindent
These first two immediately follow from the fact that $\mathcal{U} \leq \mathcal{V}_i \leq \mathcal{V}$ for all $i$, which is immediate from the definition of the $\mathcal{V}_i$.
\medskip
\item We show no element of $\mathcal{L}$ is above $deg_s(\mathcal{V})$. By discarding meets we see it suffices to show \[\bigvee_{i \in I} \mathcal{V}_i \ngeq \mathcal{V}. \] If not, for some Turing functional $\Phi_0$ we have \[ \Phi_0: \bigvee_{i \in I} \mathcal{V}_i \to \mathcal{V}. \]

\noindent
By the distributive laws we see \[ \bigvee_{i \in I} \mathcal{V}_i  \equiv (\mathcal{U} \vee \bigvee_{i \in I} \mathcal{S}_i) \wedge \mathcal{V}. \]

\noindent
Combining these last two, there is a Turing functional $\Phi_1$ so that:  \[ \Phi_1: (\mathcal{U} \vee \bigvee_{i \in I} \mathcal{S}_i) \wedge \mathcal{V} \to \mathcal{V}. \]

\noindent
Then, for each $f \in (\mathcal{U} \vee \bigvee_{i \in I} \mathcal{S}_i)$, $\Phi_1^{0 \cat f} \in \mathcal{V}$. Hence there is a Turing functional $\Phi_2$ so that $\Phi_2^f \in \mathcal{V}$ for all $f \in (\mathcal{U} \vee \bigvee_{i \in I} \mathcal{S}_i)$. So $(\mathcal{U} \vee \bigvee_{i \in I} \mathcal{S}_i) \geq \mathcal{V}$, contradicting $N_I$. 

\medskip
\item We show $\mathcal{L}$ is free. Binns and Simpson ~\cite{BiSi} pointed out that by a Lemma from Lattice Theory (for instance see ~\cite[Theorem II.2.3]{Gr}), we need only show that if \[ \bigvee_{i \in I} \mathcal{V}_i \geq \bigwedge_{j \in J} \mathcal{V}_j \] then $I \cap J \neq \emptyset.$

\noindent
Substituting and then expanding both sides by the distributive laws we have \[ (\mathcal{U} \vee \bigvee_{i \in I} \mathcal{S}_i) \wedge \mathcal{V} \geq (\mathcal{U} \vee \bigwedge_{j \in J} \mathcal{S}_j) \wedge \mathcal{V}. \]

\noindent
Next, by discarding a meet, we have \[ (\mathcal{U} \vee \bigvee_{i \in I} \mathcal{S}_i)  \geq (\mathcal{U} \vee \bigwedge_{j \in J} \mathcal{S}_j) \wedge \mathcal{V}. \] By $P_{I,J}$, we must have $I \cap J \neq \emptyset.$

\medskip
\item Finally, we must show no element of $\mathcal{L}$ is below $\deg_s(\mathcal{U})$. By discarding joins we see it suffices to show \[ \bigwedge_{j \in J} \mathcal{V}_j \nleq \mathcal{U}. \] This follows from the proof of (4) and the fact that \[ \mathcal{U} \leq \bigvee_{i \in I} \mathcal{V}_i \] for $I \cap J = \emptyset.$

\end{enumerate}

\section{Further Definitions and a Key Lemma}

It is a basic fact that if $\mathcal{U}$ is a  $\Pi_1^0$ subset of $\omega^\omega(2^\omega)$, then there is a tree $T \subset \omega^{<\omega}(2^{<\omega})$ such that $\pp = [T]$, where $[T]$ is the set of infinite paths through $T$. Sequences of such trees will be our computable approximations of $\Pi_1^0$ classes. 

We will want to consider sets of strings of fixed length in these trees, and so we have the following notation.

\begin{defn} If $T$ is a tree, let $T^s = \lbrace \sigma \in T$ : $\mid \sigma | = s \rbrace$.

\end{defn} 

Although every $\Pi_1^0$ subset $\mathcal{U}$ is the set of paths through some computable tree $T$, it might happen that $T$ contains non-extendible nodes, i.e.\ nodes that are not initial segments of any element of $\pp$. A lemma will help overcome this obstacle, but we make a definition first.

\begin{defn} A sequence of sets $\lbrace C_i \rbrace_{i \in \omega}$ is \emph{nested} if for each $j \leq k$, $C_j \supseteq C_k$. 

\end{defn}

\begin{lem} If $\mathcal{U}$ is a $\Pi_1^0$ subset of $\omega^\omega(2^\omega)$, there is a nested sequence of uniformly computable trees $\lbrace T_{\mathcal{U},s} \rbrace_{s \in \omega}$ so that each $T_{\mathcal{U},s}$ is a subset of $\omega^{<\omega} (2^{<\omega})$, $T_{\mathcal{U}} = \bigcap T_{\mathcal{U},s}$ contains only extendible nodes, and $[T_\mathcal{U}] = \pp$. 

\label{lem:apprxmtn}
\end{lem}

\begin{proof} Let $\hat{T}_{\mathcal{U}}$ be a computable tree such that $[\hat{T}_{\mathcal{U}}] = \pp$. 
\newline
Set \[ T_{\mathcal{U},s} = \lbrace \sigma \in \hat{T}_{\mathcal{U}} : \exists \tau \in \hat{T}_{\mathcal{U}}^s( \tau \subset \sigma \mbox{ or } \sigma \subset \tau) \rbrace. \] Set \[ T_{\mathcal{U}} = \bigcap T_{\mathcal{U},s}. \] It is straightforward to check $T_{\mathcal{U}}$ contains only extendible nodes and $[T_{\mathcal{U}}] = \pp$. We will call the sequence $\lbrace T_{\mathcal{U},s} \rbrace$ \emph{the canonical approximation of} $\pp$ \emph{with respect to} $\hat{T}_{\mathcal{U}}$. Usually we will take the tree $\hat{T}_{\mathcal{U}}$ for granted and simply speak of \emph{the canonical approximation}. Note that while each $T_{\mathcal{U},s}$ is computable, $T_{\mathcal{U}}$ may not be computable.

\end{proof}

If we are given two $\Pi_1^0$ subsets of $\omega^\omega$, and canonical approximations to each, we may want to have a canonical approximation for the join or meet of these subsets, in terms of our given approximations. We may even want to do this with arbitrary finite combinations of $\Pi_1^0$ subsets of $\omega^\omega$.

\begin{defn} If $\pp$ and $\pq$ are $\Pi_1^0$ subsets of $\omega^\omega$ with canonical approximations  $\lbrace T_{\pp,s} \rbrace_{s \in \omega}$ and  $\lbrace T_{\pq,s} \rbrace_{s \in \omega}$, respectively, define the \emph{canonical approximations}  $\lbrace T_{(\pp \vee \pq),s} \rbrace_{s \in \omega}$ and $\lbrace T_{(\pp \wedge \pq),s} \rbrace_{s \in \omega}$  as follows:

$\sigma \in T_{(\pp \vee \pq),s}$ iff the string formed by taking $\sigma$'s values on even numbers is in $T_{\pp,s}$ and the string formed by taking $\sigma$'s values on odd numbers is in $T_{\pq,s}$. 

$\sigma \in T_{(\pp \wedge \pq),s}$ iff there is $\tau \in T_{\pp,s}$ so that $0 \cat \tau = \sigma$ or there is $\tau \in T_{\pq,s}$ so that $1 \cat \tau = \sigma$. 

The definition can be extended inductively to give a canonical approximation of any 
\newline
\noindent
$\Pi_1^0$ subset of $\omega^\omega$ built up out of finitely many joins and meets of canonically approximated $\Pi_1^0$ subsets of $\omega^\omega$.

\end{defn}

In a straightforward way the following Lemma follows from the previous Definition.

\begin{lem} $\pp \vee \pq = [\bigcap_{s \in \omega} T_{(\pp \vee \pq),s}]$ and  $\pp \wedge \pq = [\bigcap_{s \in \omega} T_{(\pp \wedge \pq),s}].$
\end{lem} 

By induction with this lemma and the previous definition one can show that the canonical approximation of a $\Pi_1^0$ subset of $\omega^\omega$ constructed by finitely many joints and meets of $\Pi_1^0$ subsets of $\omega^\omega$ is a nested uniformly computable sequence such that the class of paths through its intersection is the intended subset.

\medskip
As we build separating classes, which are $\Pi_1^0$ subsets of $2^{\omega}$, we will want stagewise canonical approximations of them as well.

\begin{defn} Suppose $A = \bigcup_{s \in \omega} A_s$ and $ B= \bigcup_{s \in \omega} B_s$ are disjoint, c.e., and constructed in stages (i.e.\ $A_s \subseteq A_{s+1}$ and $B_s \subseteq B_{s+1}$ for all $s$). Let $\mathcal{S} = S(A, B)$. Then we define as follows a computable tree $\hat{T}_{\mathcal{S}}$  with respect to which the canonical approximation $\lbrace T_{\mathcal{S},s} \rbrace_{s \in \omega}$ for $\mathcal{S}$ is to be taken via the method in the proof of Lemma ~\ref{lem:apprxmtn}.
\newline
\[ \hat{T}_\mathcal{S} = \lbrace \sigma : \sigma(n) = 1 \Rightarrow n \notin B_{\mid \sigma |}, \sigma(n) = 0 \Rightarrow n \notin A_{\mid \sigma |} \rbrace. \]

\label{defn:S(A,B)_apprxmtn}
\end{defn}

To carry out our Sacks Preservation and Coding Strategies, we will need a way to measure incremental progress toward a final result we want to avoid.
\begin{defn}(Cenzer-Hinman ~\cite[Definition 18]{CeHi}: Length of Agreement) If $\Psi$ is a Turing functional and $\pp$ and $\pq$ are $\Pi_1^0$ classes, define:
\[ \ell ^{\pp,\pq} (\Psi, s) = \mu y[ \exists \sigma \in T^s_{\pp,s} (\Psi^{\sigma}_s \upharpoonright(y+1) \notin T_{\pq,s})]. \] If $\Psi^{\sigma}_s \upharpoonright(y+1)$ is undefined we say it is not in $T_{\pq,s}$.  
\label{dfn:length}
\end{defn}

Note that $\ell^{\pp,\pq}(\Psi,s) \geq n$ iff $\forall \sigma \in T^s_{\pp,s}(\Psi^{\sigma}_s \res n \in T_{\pq,s})$.

When it is obvious which $\pp$ and $\pq$ are under consideration, the superscripts on $\ell$ are sometimes dropped. Similarly $\Psi$ will be dropped from the argument when it is obvious which functional is under consideration.

To make the proof easier to read, we will also be interested in stages at which the length of agreement becomes greater than it ever has been, and so we have the following definitions.

\begin{defn} If $\Psi$ is a Turing functional and $\pp$ and $\pq$ are $\Pi_1^0$ classes, define \[ \overline{\ell}^{\pp, \pq}(\Psi,s) = max\lbrace \ell^{\pp, \pq}(\Psi, t) : t \leq s \rbrace. \] Or more simply, think: \[ \overline{\ell}(s) = max\lbrace \ell(t) : t \leq s \rbrace. \]

\label{dfn:ellbar}
\end{defn}

Note that $\lim \sup_s \ell(s) = \lim_s \overline{\ell}(s)$.

\begin{defn}If $\Psi$ is a Turing functional and $\pp$ and $\pq$ are $\Pi_1^0$ classes, $s + 1$ is an \emph{expansionary stage( for $\pp, \pq$, and $\Psi$)} if $\overline{\ell}^{\pp, \pq}(\Psi, s + 1) > \overline{\ell}^{\pp, \pq}(\Psi,s)$.
\end{defn}

We will act for the sake of a requirement only at an expansionary stage for the the relevant length of agreement function. This way, we easily see that if there is an upper bound on the length of agreement function for a requirement, it will act only finitely often. 

The following lemma confirms that our definition of length of agreement is well-behaved, and it will be one of the essential elements for the proof that our construction succeeds.

\begin{lem}If $\Psi$ is a Turing functional and $\pp$ and $\pq$ are $\Pi_1^0$ classes:
\begin{enumerate}
\item If $\Psi: \pp \rightarrow \pq$, then $\lim_s \ell^{\pp, \pq}(\Psi, s) = \infty$. 
\item If $\lim \sup_s \ell^{\pp, \pq}(\Psi,s) = \lim_s \overline{\ell}^{\pp, \pq}(\Psi, s) = \infty$, then $\Psi: \pp \rightarrow \pq$. 
\end{enumerate}

\label{lem:key}
\end{lem}

\begin{proof} (1) Suppose that $\Psi: \pp \rightarrow \pq$. We prove that for every $n \in \omega$ there is a stage $t$ such that for all $t' \geq t$, $\ell(t') \geq n$. Fix $n$. For each $\sigma \in T_{\pp,0}^{\mid \sigma |}$, define $\tau_\sigma = \Psi^\sigma_{ \mid \sigma |} \res n$ if it is defined. (If it is not defined, by convention say that $\tau_\sigma$ is not in any tree.) Define \[Bad = \lbrace \sigma \in T_{\mathcal{U},0} : \tau_\sigma \notin T_{\mathcal{V}, \mid \sigma |} \rbrace.\] Define \[Good = \lbrace \sigma \in T_{\mathcal{U},0} : \tau_\sigma \in T_{\mathcal{V}, \mid \sigma |} \rbrace.\] Immediately $Bad \cup Good = T_{\mathcal{U},0}$ and $Bad \cap Good = \emptyset$. Also note $Bad$ is closed downwards (and is therefore a tree) and $Good$ is closed upwards in $T_{\mathcal{U},0}.$ 

If $Bad$ were infinite, by compactness there would be an $X \in [Bad]$. Since $Bad \subset T_{\mathcal{U},0}$ and $[T_{\mathcal{U},0}] = \mathcal{U}$, we would have $X \in \mathcal{U}$. By hypothesis $\Psi: \pp \rightarrow \pq$, and so there would be $\sigma \subset X$ so that $\Psi^{\sigma}_p \res n \in T_{\mathcal{V}}$ for some stage $p$. Letting $q = max \lbrace \mid \sigma |, p \rbrace$ we have $\Psi^{X \res q}_q \res n \in T_{\mathcal{V}, q},$ whence $X \res q \in Good$, a contradiction. 

Since $Bad$ is finite, there is stage $t$ so that if $t' \geq t$ and $\sigma \in T_{\mathcal{U},t'}^{t'}$, then  $\sigma \in Good$. This means that for all $t' \geq t,$ $\ell(t) \geq n$, as desired.

(2) Suppose $\lim \sup_s \ell(s) = \infty$. Given $X \in \pp$ we want to show $Y = \Psi(X) \in \pq$. We do this by showing $Y \res n \in T_{\mathcal{V}}$ for all $n$. Since $\lim \sup_s \ell(s) = \infty$ we can find a stage $p$ so that $\ell(p) \geq n$. Then $\Psi^{X \res p}_p \res n$ is defined and  in $T_{\mathcal{V},p}$. For all stages $p' \geq p$,  $\Psi^{X \res p'}_{p'} \res n = \Psi^{X \res p}_p \res n$. Set $\tau = \Psi^{X \res p}_p \res n$. Note that $Y \res n = \tau$.

If $Y\res n = \tau \notin T_{\mathcal{V}}$, then at some stage $q$ we have $\tau \notin T_{\mathcal{V}, q}$. Because $\lbrace T_{\pq,s} \rbrace_{s \in \omega}$ is nested, for all  $q' \geq q$, $\tau \notin T_{\mathcal{V}, q'}$. Let $r \geq$ max $\lbrace q,p \rbrace$. (Actually, it is necessary that $q > p$.) For all $r' \geq r$,    $\Psi^{X \res r'}_{r'} \res n = \Psi^{X \res p}_p \res n = \tau \notin T_{\mathcal{V}, r'}$. Then $\ell(r') < n$ for all $r' \geq r$, contradicting that $\lim \sup_s \ell(s) = \infty$.

Hence $Y\res n = \Psi^{X} \res n = \Psi^{X \res p}_p \res n = \tau \in T_\mathcal{V},$ as desired.

\end{proof}

Binns ~\cite[Lemma 6]{Bi} proved something similar to part~(1) of our Lemma here. Part~(1) is actually slightly stronger than we will need for the construction. In place of part~(1) the weaker condition `If $\Psi: \pp \rightarrow \pq$ then $\lim \sup_s \ell(s) = \infty.$' would be sufficient. This weaker form of part~(1) follows immediately from Definition 18 and Proposition 19 in the proof of Cenzer and Hinman ~\cite{CeHi}. Part~(2) is actually not true for one of the length of agreement functions used in the construction of Cenzer and Hinman, and this causes the infinite injury that forces the creation of a hatted length of agreement function.

\section{Placing Restraints on $\Pi_1^0$ Classes}
\label{sec:restraints}

Suppose we are building a $\Pi_1^0$ class $\mathcal{S} = S(A,B)$ by building $A$ and $B$ as c.e.\ subsets of $\omega$ such that $A \cap B = \emptyset$ and certain other requirements are met. As indicated in Definition ~\ref{defn:S(A,B)_apprxmtn}, there is a canonical approximation $\lbrace T_{\mathcal{S}, s} \rbrace_{s \in \omega}$ so that $\mathcal{S} = [\bigcap_{s \in \omega} T_{\mathcal{S},s}].$ 

For a strategy at stage $t$ to restrain $\mathcal{S}$ up to level $m$ means to attempt to ensure that for all stages $t' \geq t$, if $\mid \tau| \leq m$ and $\tau \in T_{\mathcal{S},t}$ then $\tau \in T_{\mathcal{S}, t'}$. 

The strategy will force lower priority strategies to comply with this request. The restraint may fail (`be injured') if a higher priority strategy makes an enumeration that violates this request. In particular, since $\mathcal{S} = S(A,B)$, if $\tau \in T_{\mathcal{S},t}$ but $\tau \notin T_{\mathcal{S},t'}$, it must be that for some $x < \mid\tau|$, either

\begin{enumerate}
\item $\tau(x) = 1$, but $ x \in B_{t'} \setminus B_t$ or 
\item $\tau(x) = 0$, but $x \in A_{t'} \setminus A_t$. 
\end{enumerate}

In other words, restraining $\mathcal{S}$ up to level $m$ at stage $t$ amounts to an attempt to prevent the following situations: $ x \in B_{t'} \setminus B_t$ or $x \in A_{t'} \setminus A_t$ for some $t' \geq t$ and some $x < m$. This amounts to restraining $A_t \res m$ and $B_t \res m$ in the sense standard for c.e.\ priority arguments. 

If we are told that we need not worry about anything beyond ensuring that at expansionary stages $t$, $A$ is protected from injury by lower priority requirements up to some specified level $m$ (i.e. $A_t \res m = A \res m$, if there is no injury by higher priority requirements), then the work is even simpler. For if $A_t(x) = 1$ for some $x < m$, then there is no way we will enumerate $x \in B_{t'}$ at a later stage $t' \geq t$, because we insist $A \cap B = \emptyset$. Therefore, case (1) from above will never be a concern, and we need not restrain $B$ at all. This simplified method for placing restraints is the one we will use in our construction.

\section{Negative Requirements and Strategies}
\label{sec:negative}

To satisfy $N_I$, we will have for each Turing functional $\Phi$ the requirement 
\[ N_{I, \Phi} \hspace{0.4cm} \Phi: \mathcal{U} \vee \bigvee_{i \in I} \mathcal{S}_i \nrightarrow \mathcal{V}. \]

\bigskip
A simplified version of the Sacks Preservation Strategy is used for satisfying the negative requirements. At expansionary stages we will restrain each $A_i$ at least up to its use in the relevant computations. 

\begin{enumerate}

\item Wait for an expansionary stage $s$. 
\item For each $i \in I$ restrain $A_i$ up to its maximum use in all computations used in calculating $\ell ^{(\mathcal{U} \vee \bigvee_{i \in I}  \mathcal{S}_i),\pq} (\Phi, s).$ We may as well take all the restraints to be $s+1$. Also, initialize all $P_{I, J, \Psi}$ strategies of lower priority than $N_{I, \Phi}$ with markers $m_{\sigma, j} < s+1$. 

\item Go back to step 1 and wait for another expansionary stage.

\item $N_{I, \Phi}$ is injured if a higher priority positive requirement performs an enumeration that violates an $A_i$-restraint. No action is taken.

\end{enumerate}

\bigskip
\textbf{Current Outcome} at stage $s$ is $\overline{\ell}(s)$. 

\textbf{Final Outcome} is $\lim_s \overline{\ell} (s) = \lim \sup_s \ell (s)$.  

\bigskip
\noindent
\underline{Verification}: that $N_{I, \Phi}$ is satisfied and acts only finitely often.

Suppose $N_{I,\Phi}$ is not satisfied. Then $\Phi: \pp \vee \bigvee_{i \in I} \mathcal{S}_i \rightarrow \pq$. By Lemma 4.2(1) $\lim_s \ell (s) = \infty$. Hence $\lim_s \overline{\ell}(s) = \lim \sup_s \ell(s) = \infty$, and there were infinitely many expansionary stages. 

By induction assume that higher priority requirements act only finitely often. Then after some stage $q$ no $A_i$-restraints are ever injured. We derive a contradiction to the theorem's hypothesis by giving a uniform procedure to calculate a $Y \in \pq$ given $X \in \pp$. 

Fix $n$. We describe how to calculate $Y \res n$.  Look for the first expansionary stage $t \geq q$ so that $\ell^{(\pp \vee \bigvee_{i \in I} \mathcal{S}_i),\mathcal{V}}(\Phi,t) \geq n$. Such a $t$ exists because there were infinitely many expansionary stages. Define \mbox{$\sigma_t = (X \oplus \bigoplus_{i \in I} A_{i,t}) \res t.$} Note that $\sigma_t \in T^t_{(\mathcal{U} \vee \bigvee_{i \in I} \mathcal{S}_{i}), t}$. Therefore, by the definition of the length of agreement function, $\Phi^{\sigma_t}_t \res n \in T_{\mathcal{V},t}$. 

Set $Y \res n = \Phi_t^{\sigma_t} \res n$. At stage $t$ step (2) of the strategy directs us to restrain each $A_i$ up to level $t+1$, which is greater than the use of each $A_i$ in the computation showing $\Phi^{\sigma_t}_t \res n \in T_{\mathcal{V},t}$. By choice of $q$ the restraints up to $t+1$ on each $A_i$ will never be violated. Therefore,  
\newline
$\Phi_t^{\sigma_t} \res n = \Phi^{X \oplus \bigoplus_{i \in I} A_i} \res n$. Furthermore, $\Phi^{X \oplus \bigoplus_{i \in I} A_i} \res n \in T_{\pq}$, because $X \oplus \bigoplus_{i \in I} A_i \in \pp \vee \bigvee_{i \in I} \mathcal{S}_i$ and $\Phi: \pp \vee \bigvee_{i \in I} \mathcal{S}_i \rightarrow \pq$. Hence $Y \res n \in T_{\pq}$ for each $n$, and so $Y \in \pq$, as desired. 

This contradiction shows the requirement is in fact satisfied. By the contrapositive of Lemma ~\ref{lem:key}(2), there are only finitely many expansionary stages, and so $N_{I, \Phi}$ acts only finitely often to impose restraints on each $A_i$. 

\subsection{Contrast with Negative Requirements in the Proof of the Density of $\mathcal{P}_s$.}
\label{subsec:CenHin}

Our negative requirements, met by Sacks Preservation Strategies, were:

\[ N_{I, \Phi} \hspace{0.4cm} \Phi: \mathcal{U} \vee \bigvee_{i \in I} \mathcal{S}_i \nrightarrow \mathcal{V}. \]

Cenzer and Hinman had negative requirements very similar to these ~\cite[p. 590]{CeHi}. In our notation, they began with $\mathcal{U} < \mathcal{V}$ and were building $\mathcal{S} = S(A,B)$, where $A$ and $B$ are disjoint c.e\ subsets of $\omega$ built by the construction. The requirements ensured that 
\newline
$\mathcal{U} < (\mathcal{U} \vee \mathcal{S}) \wedge \mathcal{V} < \mathcal{V}.$ For each Turing functional $\Phi$, their negative requirement was (in our notation),

\[ N_{\Phi} \hspace{0.4cm} \Phi: \mathcal{U} \vee \mathcal{S} \nrightarrow \mathcal{V}. \] 

These are the same as our negative requirements, except we have the join of finitely many $\mathcal{S}_i$ in place of one $\mathcal{S}$. Cenzer and Hinman sought to simplify, by replacing each $N_{\Phi}$ with a requirement:

\[N^A_{\Phi} \hspace{0.4cm} \neg \forall X \in \mathcal{U}(\Phi^{X \oplus A} \in \mathcal{V}). \]

Because $A \in S(A, B) = \mathcal{S}$, the satisfaction of $N^A_{\Phi}$ guarantees the satisfaction of $N_{\Phi}$. In a way, $N_{\Phi}^A$ is a simpler requirement than $N_{\Phi}$: there is less to keep track of. However, in general $\lbrace A \rbrace$ is not a $\Pi_1^0$ class. Therefore, the length of agreement function $\ell$ from Definition ~\ref{dfn:length} cannot be directly adapted to work for $N_{\Phi}^A$.

This approach of Cenzer and Hinman would correspond for us to requirements of the form \[N_{I, \Phi}^{A_i} \hspace{0.4cm} \neg \forall X \in \mathcal{U}(\Phi^{X \oplus \bigoplus_{i \in I} A_i} \in \mathcal{V}). \] Again, because $A_i \in S(A_i, B_i) = \mathcal{S}_i$ for each $i$, the satisfaction of $N_{I, \Phi}^{A_i}$ guarantees the satisfaction of $N_{I, \Phi}$. Of course, the length of agreement function $\ell$ cannot be directly adapted for $N^{A_i}_{I,\Phi}$ either. 

Because the original length of agreement function was not directly adaptable to the demands of the new requirement, Cenzer and Hinman defined another length of agreement function ~\cite[pp. 594-595]{CeHi}. This length of agreement function works directly with the c.e.\ set $A$:

\[ \ell^{(\mathcal{U} \times A), \mathcal{V}}  (\Phi,s) = \mu y[(\exists \sigma \in T^s_{\mathcal{U},s}) \hspace{0.18cm}\hat{\Phi}^{\sigma, A_s}_s \res (y + 1) \notin T_{\mathcal{V},s}], \]
\noindent
where $\hat{\Phi}$ is defined via the hat trick adapted for $\Pi_1^0$-classes. The hat trick is needed to handle the infinite injury that accompanies this new length of agreement function. For it can happen that $\lim \sup_s \ell^{(\mathcal{U} \times A), \mathcal{V}}  (\Phi,s) = \infty$, although $N_{\Phi}^A$ is satisfied. This is precisely where Lemma ~\ref{lem:key}(2) fails, as mentioned after the proof of that Lemma.

See Soare \cite[Chapter 8]{So} for an explanation of the original hat trick as used to combat infinite injury in a proof of the Sacks Density Theorem.  It seems the hat trick is necessary for a direct proof of the Sacks Density Theorem by priority argument.

This strategy of Cenzer and Hinman using the hat trick would also work in our case. With only straightforward extensions of definitions, requirements of the form $N^{A_i}_{I, \Phi}$ can be satisfied for us without need of any further work.

However, it is interesting to note that infinite injury and the hat trick machinery can be avoided as we do in this paper, by working with the original requirements, namely $N_{I, \Phi}$ in our case, and $N_{\Phi}$ for Cenzer and Hinman. The construction and verification used in our negative strategy goes through without any extra work when applied to the $N_{\Phi}$-requirements of Cenzer and Hinman.

\section{Positive Requirements and Strategies}
To satisfy $P_{I, J}$ we have for each Turing functional $\Psi$ the requirement, 

\[ P_{I,J, \Psi} \hspace{0.4cm} \Psi: \mathcal{U} \vee \bigvee_{i \in I} \mathcal{S}_i \nrightarrow (\mathcal{U} \vee \bigwedge_{j \in J} \mathcal{S}_j) \wedge \mathcal{V} \] for $I, J$ finite, $I \cap J = \emptyset$. 

If $P_{I, J, \Psi}$ fails we ensure there are Turing functionals
\[ \Gamma_j: \mathcal{S}_j \to \mathcal{V} \] for each $j \in J$ and 
\[ \Delta_i: \textbf{0} \to \mathcal{S}_i \] for each $i \in I$. 

\bigskip
Then we can put $\Psi, \Gamma_j$, and $\Delta_i$ together to show $\mathcal{U} \geq \mathcal{V},$ contradicting the hypothesis of the theorem. Here's how: given $X \in \pp$ we can effectively produce $X_i \in \mathcal{S}_i$ for each $i$ via the finitely many $\Delta_i$. Then $\Psi (X \oplus \bigoplus_{i \in I} X_i)$ produces either $1 \cat Z$ where $Z \in \pq$ or $0 \cat Z$ where \[ Z \in \pp \vee \bigwedge_{j \in J} \mathcal{S}_j. \] If we see $1 \cat Z$, then $Z \in \pq$, and we are done. If we see $0 \cat Z$, then given the finite piece of information $|J|$, we can extract from $Z$ a $W \in \mathcal{S}_j$ for some $j \in J$, and we will know which $j$ it is. Then $\Gamma_j (W)$ gives $Y \in \pq$ as desired.

\bigskip
To create the Turing functionals $\Gamma_j$ we use a Sacks Coding Strategy. To create the Turing functionals $\Delta_i$, in a way similar to a Sacks Preservation Strategy we make each $A_i$ computable via restraints set at expansionary stages.  We proceed as follows: 

\bigskip

\begin{enumerate}

\item Set $n = 0$.

\item Wait for an expansionary stage $s$ at which $\ell^{(\mathcal{U} \vee \bigvee_{i \in I} \mathcal{S}_i), ((\mathcal{U} \vee \bigwedge_{j \in J} \mathcal{S}_j) \wedge \mathcal{V})}(\Psi, s) \geq n$. 

\item For each $i \in I$ restrain $A_i$ up to level $s + 1$. Note that $s + 1$ will be greater than the length of agreement at stage $s$. Also, initialize all $P$-strategies of lower priority with markers $m_{\sigma,j} < s + 1$. 

\item For each $j \in J$ and each $\sigma \in 2^{<\omega}$ with $\mid \sigma | = n$ choose large markers  $m_{\sigma, j}$. Here ``large'' means greater than any restraints on $A_j$ and greater than any other markers for $\mathcal{S}_j$ (whether or not at this stage they have been enumerated into $A_j$ or $B_j$) established for the sake of this or any other $P$-strategy. Increment $n$ and go to step 2 for this new $n$. Meanwhile, for the old $n$ go to the next step.

\item Wait for another expansionary stage $t > s$ at which either
\newline
(i) $\sigma \cat 0 \notin T_{\mathcal{V}, t}$ or 
\newline
(ii) $\sigma \cat 1 \notin T_{\mathcal{V}, t}$.

\item If (i) happens for $\sigma$ at stage $t$, place $m_{\sigma, j}$ into $A_j$ for each $j \in J$, and do not act again for this $m_{\sigma, j}$. 

\item If (ii) happens for $\sigma$ at stage $t$, place $m_{\sigma, j}$ into $B_j$ for each $j \in J$ and do not act again for this $m_{\sigma, j}$.

\item $P_{I,J,\Psi}$ is injured if a higher priority $P$-strategy makes an enumeration that violates an $A_i$-restraint. No action is taken.

\end{enumerate}

\bigskip
As directed in the strategies, a $P_{I, J, \Psi}$ strategy is initialized if a higher priority $N$- or \mbox{$P$-strategy} makes an $A_j$-restraint greater than a current marker $m_{\sigma, j}$ for $P_{I, J, \Psi}$. Note that if some requirement $N_{I', \Phi}$ or $P_{I', J', \Psi'}$ initializes $P_{I, J, \Psi}$, it must be that $I' \cap J \neq \emptyset$.
\newline
\indent
If $P_{I, J, \Psi}$ is initialized, all current markers $m_{\sigma, j}$ are discarded and we return to step 1, setting $n = 0$. 
\newline
\indent
The $P_{I,J,\Psi}$ outcomes are the same as the outcomes for negative strategies:

\bigskip
\noindent
\textbf{Current Outcome} at stage $s$ is $\overline{\ell} (s).$
\newline
\textbf{Final Outcome} is $\lim_s \overline{\ell} (s) = \lim \sup_s \ell (s)$.  

\bigskip
\noindent
\underline{Verification}: that $P_{I,J,\Psi}$ is satisfied and acts only finitely often. 

First note that the action for the $\lbrace \mathcal{S}_i \rbrace_{i \in I}$ and the $\lbrace \mathcal{S}_j \rbrace_{j \in J}$ never conflict because $I \cap J = \emptyset$. 

By induction assume all higher priority requirements have finite outcomes, so that there is a stage $q$ after which $P_{I, J, \Psi}$ is never initialized and is never injured.

\bigskip
If $P_{I, J, \Psi}$ were not satisfied, we would have by Lemma 2.2(1) that $\lim_s \overline{\ell} (s) = \lim \sup_s \ell (s) = \infty$. In this case, we need to define the functionals $\Gamma_j$ and $\Delta_i$ which will yield a contradiction. We begin by giving $\Gamma_j: \mathcal{S}_j \to \pq$ for each $j \in J$. 

Given $X \in \mathcal{S}_j$ and inductively assuming we have $\sigma_n \in T_{\mathcal{V}}$ with $\mid \sigma_n| = n$, we say how to determine if $\sigma_n \cat 0 \in T_{\mathcal{V}}$ or $\sigma_n \cat 1 \in T_{\mathcal{V}}$. First determine the value, if any, of $m_{\sigma_n, j}$ after the last initialization sometime before stage $q$. If $m_{\sigma_n, j}$ was not defined by stage $q$ simply find its value once it is defined. Note that because $\lim_s \overline{\ell} (s) = \infty$, $m_{\sigma_n, j}$ must eventually be defined. Further, $m_{\sigma_n, j}$ can never change after stage $q$. 

\bigskip
\underline{Case 1} If $X(m_{\sigma_n, j}) = 1$, then $m_{\sigma_n, j} \notin B_j$ because $X \in S(A_j, B_j)$. Thus if (ii) happened from step 4 in the strategy, it happened after (i). Since $T_{\mathcal{V}}$ contains only extendible nodes and $\sigma_n \in T_{\mathcal{V}}$, we cannot have both (i) and (ii) happening. Thus (ii) never happened. Hence $\sigma_n \cat 1 \in T_{\mathcal{V}}$. Set $\sigma_{n+1} = \sigma_n \cat 1.$

\bigskip
\underline{Case 2} If $X(m_{\sigma_n, j}) = 0$, then $m_{\sigma_n, j} \notin A_j$ because $X \in S(A_j, B_j)$. Thus if (i) happened from step 4 in the strategy, it happened after (ii). Since $T_{\mathcal{V}}$ contains only extendible nodes and $\sigma_n \in T_{\mathcal{V}}$, we cannot have both (i) and (ii) happening. Thus (i) never happened. Hence $\sigma_n \cat 0 \in T_{\mathcal{V}}$. Set $\sigma_{n+1} = \sigma_n \cat 0$. 

\bigskip
Defining $Y = \bigcup_n \sigma_n$, we see that because each of its initial segments is in $T_{\mathcal{V}}$, we have $Y \in \pq$ as desired. 

\bigskip
We now give $\Delta_i: 0 \to \mathcal{S}_i$ for each $i \in I$. We do this by calculating $A_i \in \mathcal{S}_i = S(A_i, B_i)$. To calculate $A_i \res x$ wait for an expansionary stage $t > q$ so that $\ell(t) = \overline{\ell}(t) \geq x.$ At stage $t$ in step (3) of the strategy we will restrain $A_i$ up to level $t+1$, which must be greater than $x$. Recall that after stage $q$,  $P_{I, J, \Psi}$ is never injured. So $A_i \res x$ is never injured after stage $t$.  Hence $A_{i, t} \res x = A_i \res x$. 

\bigskip
$\Gamma_j$  and $\Delta_i$ for each $j \in J$ and $i \in I$ together with \[\Psi: \mathcal{U} \vee \bigvee_{i \in I} \mathcal{S}_i \to (\mathcal{U} \vee \bigwedge_{j \in J} \mathcal{S}_j) \wedge \mathcal{V} \]
will give a uniform reduction from $\pp$ to $\pq$ as shown above, which is a contradiction. 

\bigskip
The contradiction we have arrived at shows that $P_{I, J, \Psi}$ is satisfied. Therefore, by the contrapositive of Lemma ~\ref{lem:key}(2), $P_{I, J, \Psi}$ has finite outcome. Because it acts only at expansionary stages, $P_{I, J, \Psi}$ acts only finitely often. 

Because there will be no infinite outcomes, any priority ordering for the requirements will yield a construction in which the strategies satisfy their requirements as described above, thus proving the theorem. 

\section{Appendix I: A Note on Preservation Strategies}

To accomplish our preservation strategy as part of satisfying a positive requirement $P_{I, J, \Psi}$, we simply ensure each $A_i$ is computable, if it turns out the requirement is violated. This helps lead to a contradiction of the theorem's hypothesis. 

On the other hand, in the verification that our preservation strategy satisfies a negative requirement $N_{I, \Phi}$, we don't say each $A_i$ must be computable whenever the requirement is violated. However, it is true, so long as we take the restraints at each expansionary stage $s$ to be $s + 1$, as mentioned in step (2) of the strategy for negative requirements. The verification that each $A_i$ is computable is exactly the same as in the verification for the positive requirements.  Furthermore, the computability of each $A_i$ would be enough to reach the contradiction of the theorem's hypothesis in the verification of the negative requirement. 

We leave the verification of the strategy for the negative requirements as is, because it is then as close as possible to the verification of the negative strategy of Cenzer and Hinman, with the exception of the modification that eliminates infinite injury. This way, it is easier to isolate the exact cause of infinite injury for Cenzer and Hinman.

\section{Appendix II: More on the Length of Agreement Function.}
\label{sec-comparison}

This section gives a more precise account of exactly how the length of agreement function behaves in this paper. It also applies to the construction of Cenzer and Hinman once the modification we describe in Section ~\ref{subsec:CenHin} is made. This section is not necessary for the proof of the main result of this paper.

We begin by noticing that if $\lim \sup_s \ell(s) = \infty$, then by two applications of Lemma ~\ref{lem:key}, $\lim_s \ell(s) = \infty$. Then $\lim \inf_s \ell(s) = \infty$. Hence we have the following Proposition.

\begin{prop} For a length of agreement function $\ell$ as defined in Definition ~\ref{dfn:length}, the following are equivalent.

\begin{enumerate}

\item  $\lim \sup_s$ $\ell(s) = \infty$. 
\item $\lim \inf_s$ $\ell(s) = \infty$. 
\item $\lim_s$ $\ell(s) = \infty$. 

\end{enumerate}
In particular, it is never the case that $\lim \inf_s$ $\ell(s)$ is finite but $\lim \sup_s$ $\ell(s) = \infty$. 
\label{prop:ellinf}
\end{prop}

Proposition ~\ref{prop:ellinf} leaves open the possibility that $\lim \sup_s \ell(s)$ and $\lim \inf_s \ell(s)$ are both finite but $\lim \inf_s \ell(s) < \lim \sup_s \ell(s)$. However, $\Pi_1^0$ classes are nice enough that this never happens. 

\begin{prop} For the length of agreement function $\ell$ defined in Definition ~\ref{dfn:length}, 
\newline
$\lim \inf_s$~$\ell(s) = \lim \sup_s$~$\ell(s) = \lim_s$ $\ell(s)$. These equal limits may be finite or infinite. 

\label{prop:nicelength}
\end{prop}

For the proof of this Proposition, we need a further Lemma. As motivation suppose that at some first stage $t$, $\ell(t) \geq n$. We want to analyze what could later cause the length of agreement to drop below $n$ after $t$, and also what could bring the length of agreement back up to at least $n$, once it drops below $n$. 

\begin{lem} Let $\mathcal{M}, \mathcal{N}$ be $\pie$ classes. Let $\Psi$ be a Turing functional. Fix $n \in \omega$. There is a stage $r$ such that either 

\begin{enumerate}
\item for all $r' \geq r$, $\ell^{\mm, \nn}(\Psi, r') \geq n$ or
\item for all $r' \geq r$, $\ell^{\mm, \nn}(\Psi, r') < n$. 

\end{enumerate}
\label{lem:essenceoflength}
\end{lem}

Note that Lemma ~\ref{lem:essenceoflength} immediately implies Proposition ~\ref{prop:nicelength}.

\begin{proof}

In this proof, if $\sigma \in 2^{<\omega}$ and $x \leq \mid\sigma|$ then $\sigma_x = \sigma \upharpoonright x$.

We begin with an easy fact: if $\lbrace C_i \rbrace_{i \in \omega}$ is a nested sequence and for each $i \in \omega$, $C_i$ is finite, then there is $j \in \omega$ so that for all $j' \geq j$, $C_{j'} = C_j$. 

Suppose there is a stage $t$ so that $\ell^{\mm, \nn}(\Psi, t) \geq n$. (Otherwise, for any $r$, (2) holds). Let $Output_{n,t} = \lbrace \tau \in T_{\nn, t} : \tau = \Psi^\sigma_t \upharpoonright n$ for some $\sigma \in T^t_{\mm, t} \rbrace.$ For $t' \geq t$, define $Output_{n,t'} = \lbrace \tau \in Output_{n,t}: \tau \in T_{\nn, t'} \rbrace$. By the easy fact, there is a stage $q$ so that for all $q' \geq q$, $Output_{n,q} = Output_{n,q'}$. Also by the easy fact, there is $r \geq q$ so that $T^t_{\mm,r} = T^t_{\mm,r'}$ for all $r' \geq r$. 

We claim this is the desired $r$. We verify the claim in two cases.

\medskip

Suppose $\ell^{\mathcal{M},\mathcal{N}}(\Psi,r) \geq n$. 
\newline
Let $r' \geq r$ and $\sigma \in T^{r'}_{\mathcal{M},r'}$. Because $\ell^{\mathcal{M}, \mathcal{N}}(\Psi,t) \geq n$, $\Psi^{\sigma_t}_t \upharpoonright n \in T_{\nn,t}$. Since $\Psi^{\sigma_t}_t \upharpoonright n$ converges and $r' \geq r \geq t,$ we have
$\Psi^\sigma_{r'} \res n = \Psi^{\sigma_r}_r \res n = \Psi^{\sigma_t}_t \res n \in T_{\nn,t}$.
Then $\Psi^{\sigma_r}_r \res n \in Output_{n,t}$. Since $\ell^{\mm, \nn}(\Psi,r) \geq n$, we must have $\Psi^{\sigma_r}_r \res n \in T_{\nn,r}$. So $\Psi_r^{\sigma_r} \res n \in Output_{n,r}$. Since $r \geq q$, by choice of $q$ above we have $Output_{n,r'} = Output_{n,r}$, and so $\Psi^{\sigma_r}_r \res n \in T_{\nn, r'}.$ But $\Psi^\sigma_{r'} \res n = \Psi^{\sigma_r}_r \res n$, so $\Psi^\sigma_{r'} \res n \in T_{\nn,r'}$. 
\newline
Hence $\ell^{\mm, \nn}(\Psi, r') \geq n$. 

\medskip

Suppose $\ell^{\mm, \nn}(\Psi,r) < n$. 
\newline
Then there is $\sigma \in T^r_{\mm,r}$ so that $\Psi^\sigma_r \res n \notin T_{\nn, r}$. But $\ell^{\mm, \nn}(\Psi,t) \geq n$, so we must have $\Psi^{\sigma_t}_t \res n \in T_{\nn, t}$. Then $\Psi^\sigma_r \res n$ must converge and $\Psi^\sigma_r \res n = \Psi_t^{\sigma_t} \res n.$ So we must have $\Psi_t^{\sigma_t} \res n \notin T_{\nn, r}$. Because $\sigma \in T^r_{\mm, r}$, we must have $\sigma_t \in T^t_{\mm,r}$. By the choice of $r$, we must have $\sigma_t \in T^t_{\mm,r'}$ for all $r' \geq r$. Then, by the definition of the canonical approximation (see the proof of Lemma ~\ref{lem:apprxmtn}), for all $r' \geq r$ there is $\sigma' \in T^{r'}_{\mm, r'}$ so that $\sigma' \supseteq \sigma_t$. $\Psi^{\sigma'}_{r'} \res n = \Psi_t^{\sigma_t} \res n \notin T_{\nn,r}$. Because $\lbrace T_{\nn, s} \rbrace_{s \in \omega}$ is nested, $\Psi_{r'}^{\sigma'} \res n \notin T_{\nn, r'}$. 
\newline
Hence, $\ell^{\mm, \nn}(\Psi,r') < n$, for all $r' \geq r$. 

\end{proof}

\end{document}